\documentclass[10pt,a4paper,reqno]{amsart}
\usepackage{amssymb,amsmath,amsthm}
\usepackage[T1]{fontenc}
\usepackage{graphicx}

\theoremstyle{plain}
\newtheorem{theorem}{Theorem}
\newtheorem{corollary}{Corollary}
\newtheorem{lemma}{Lemma}

\newtheorem*{theorem A}{Theorem A}
\newtheorem*{theorem B}{Theorem B}
\newtheorem*{corollary A}{Corollary A}
\newtheorem*{corollary B}{Corollary B}
\theoremstyle{definition}
\newtheorem{definition}{Definition}
\newtheorem{rem}{Remark}

\newtheorem*{proof A}{Proof of Theorem A}
\newtheorem*{proof B}{Proof of Theorem B}
\newtheorem*{proof AA}{Proof of Corollary A}

\DeclareMathOperator*{\hd}{\mathrm{dim}_\mathrm{H}}

\title{Lipschitz continuity of the Hausdorff dimension of self-affine sponges at Sierpinski sponges}
\author{Nuno Luzia}
\address{Universidade Federal do Rio de Janeiro, Instituto de Matem\'atica \\ Rio de Janeiro 21941-909, Brazil}  
\email{nuno@im.ufrj.br} 
\begin{document} 
\begin{abstract}
The Hausdorff dimension of general Sierpinski carpets, \cite{4} and \cite{20}, and the generalization on Lalley-Gatzouras carpets, \cite{10}, are today well known results, the formulas being obtain via the variational principle for the dimension. We call the multidimensional versions of these carpets Sierpinski sponges and self-affine sponges, respectively,. In this paper we show that the Hausdorff dimension of self-affine sponges, defined in $\mathrm{R}^3$, is a Lipschitz continuous function at Sierpinski sponges.
\end{abstract}
\maketitle

\section{Introduction and statements}

The dimension theory of $C^{1+\alpha}$ conformal repellers is well understood by means of the thermodynamic formalism introduced by Sinai-Ruelle-Bowen \cite{24}, \cite{23}, \cite{5} and the famous Bowen's equation \cite{6}, \cite{22}. In particular there is a unique ergodic measure of full dimension which is a Gibbs state.

The dimension theory of \emph{non-conformal} repellers is still being developed and no such general formalism exists. The computation of Hausdorff dimension of non-conformal fractals began with the fundamental works by Bedford \cite{4} and McMullen \cite{20} on the \emph{general Sierpinski carpets}, and their generalization \cite{10} on the \emph{Lalley-Gatzouras carpets}, as they are known today. See also \cite{1} for an extension of the Lalley-Gatzouras carpets. These are self-affine fractals in the plane and there is an ergodic measure of full dimension, in fact Bernoulli (such a measure is, in general, not \emph{unique} for Lalley-Gatzouras carpets \cite{2}). There are also some \emph{non-linear} versions of these results, \cite{11}, \cite{15}, \cite{17} and \cite{19}, towards a \emph{Dimension Formalism} for these kind of non-conformal repellers. In particular, in \cite{15} and \cite{19} we compute the Hausdorff dimension of \emph{Non-linear Lalley-Gatzouras carpets} which, as the name suggests, are the $C^{1+\alpha}$ non-linear versions of the Lalley-Gatzouras carpets.

Other approaches try to give a formula for the Hausdorff dimension in a \emph{generic} setting, instead of considering particular families as before. One such example is the famous \emph{Falconer's formula} \cite{8} which gives, under some hypotheses, the Hausdorff dimension of the self-affine fractal for `almost every' translation vector-parameters, but does not tell for which parameters it holds (exceptions being the self-affine fractals in the plane considered in \cite{12} and its non-linear versions \cite{16}). In fact these two approaches are quite different, for the self-affine fractals for which Falconer's formula holds there is coincidence with Hausdorff and box-counting dimensions, as for the Lalley-Gatzouras carpets these two dimensions do not coincide in general \cite{10}. More recently, Hochmann and Rapaport \cite{13} gave a formula for the Hausdorff dimension of self-affine sets in the plane satisfying an \emph{irreducibility} condition and  an \emph{exponential separation} condition.

The computation of the Hausdorff dimension of non-conformal fractals in $\mathbb{R}^d$, $d>2$ reveals to be a much more difficult task. There is a natural way of defining the $d$-dimensional versions of general Sierpinski carpets and Lalley-Gatzouras carpets which we shall call, respectively, \emph{Sierpinski sponges} and \emph{self-affine sponges}.
Recently, in a major breakthrough, \cite{7} showed that the \emph{variational principle for the dimension} does not hold, in general, within the class of \emph{Baranski sponges}  (even for $d=3$), i.e. there is not an ergodic measure of full dimension. They showed that the Hausdorff dimension in that class can be calculated via \emph{pseudo-Bernoulli measures}. On the other hand, Kenyon and Peres \cite{14} computed the Hausdorff dimension of Sierpinski sponges by proving the variational principle for the dimension and, moreover, there is a unique ergodic measure of full Hausdorff dimension, in fact Bernoulli (see also \cite{9}, and see \cite{3} for a random version). We do not know if the variational principle for the dimension holds for self-affine sponges that are close to Sierpinski sponges but we will show that the Hausdorff dimension of self-affine sponges, $d=3$, is Lipschitz continuous at Sierpinski sponges. 
Since the variational principle for the dimension holds for Sierpinski sponges, this implies that the variational principle for the dimension \emph{almost} holds for self-affine sponges close to Sierpinski sponges. We notice that the continuity of the Hausdorff dimension in the class of \emph{Baranski sponges} was proved in \cite{7}, but the class of self-affine sponges considered in this paper need not be \emph{Baranski sponges}. Also by restricting to the continuity of Hausdorff dimension at Sierpinski sponges we are able to obtain more explicit and quantitative results. 

\subsection{Sierpinski sponges}
We begin by describing the \emph{Sierpinski sponges} studied in \cite{14} (the multidimensional versions of the general Sierpinski carpets). Let $\mathbb{T}^d=\mathbb{R}^d/\mathbb{Z}^d$ be the $d$-dimensional torus and $f\colon \mathbb{T}^d\to \mathbb{T}^d$ be given by
\[
    f(x_1,.x_2,...,x_d)=(l_1 x_1, l_2 x_2,...,l_d x_d)
\]
where $l_1\ge l_2\ge...\ge l_d>1$ are integers. The grids of hyperplanes
\begin{align*}
       \{i/l_1\}\times [0,1]^{d-1}, & \,i=0,..., l_1-1\\
       [0,1]\times \{i/l_2\}\times [0,1]^{d-2}, &\, i=0,..., l_2-1\\
       &\vdots\\
       [0,1]^{d-1}\times\{i/l_d\}, & \,i=0,..., l_d-1
\end{align*}
form a set of boxes each of which is mapped by $f$ onto the entire torus (these boxes are the domains of invertibility of $f$). Now choose some of these boxes and
consider the fractal set $\Lambda$ consisting of those points that always remain in these chosen boxes when iterating $f$. Geometrically, $\Lambda$ is the limit (in 
the Hausdorff metric), or the intersection, of \emph{$n$-approximations}: the 1-approximation consists of the chosen boxes; the 2-approximation consists in replacing each box of the 1-approximation by a rescaled affine copy of the 1-approximation, resulting in more and smaller boxes;  the 3-approximation consists in replacing each box of the 2-approximation by a rescaled affine copy of the 1-approximation, and so on. We say that $\Lambda$ is a \emph{Sierpinski sponge}, and their Hausdorff dimension was computed in \cite{14} (a formula is given in next section).

\subsection{Self-affine sponges}
Now we describe the \emph{self-affine sponges} (the multidimensional versions of Lalley-Gatzouras carpets), which include as a very particular case the Sierpinski sponges.
Let $S_1, S_2, ..., S_r$ be contractions of $\mathbb{R}^d$. Then there is a unique nonempty compact set $\Lambda$ of $\mathbb{R}^d$ such that
\[
   \Lambda=\bigcup_{i=1}^{r} S_i(\Lambda) .
\]
We will refer to $\Lambda$ as the limit set of the semigroup generated by $S_1, S_2, ..., S_r$. We are going to consider sets $\Lambda$ which are limit sets of the semigroup generated by the $d$-dimensional mappings $A_{i^1i^2...i^d}$ given by
\[
 A_{i^1i^2...i^d}=\begin{pmatrix} 
  a_{i^1i^2...i^d} & 0 & 0 & \cdots & 0 \\ 0 & a_{i^1i^2...i^{d-1}} & 0 & \cdots & 0 \\ & \vdots & \\ 0 & \cdots & 0 & 0 & a_{i^1} 
  \end{pmatrix}x +  \begin{pmatrix} u_{i^1i^2...i^d} \\ u_{i^1i^2...i^{d-1}} \\ \vdots \\ u_{i^1}
  \end{pmatrix}
\]
for $(i^1,i^2,...,i^d)\in \mathcal{I}$. Here
\begin{align*}
 \mathcal{I}=\{(i^1,i^2,...,i^d):\;&1\le i^1\le m\,,\, 1\le i^2\le m_{i^1},\\
   &1\le i^3\le m_{i^1i^2}\,,...,\, 1\le i^d\le m_{i^1i^2...i^{d-1}} \}
\end{align*}
is a finite index set, and $0<a_{i^1...i^{k}}<1\,,\,k=1,...,d$ satisfy
\[
   a_{i^1...i^ki^{k+1}}\le a_{i^1...i^k}\,.
\]
Also, for each $(i^1,...,i^d)\in\mathcal{I}$ and $k\in\{1,...,d\}$,
\[
   \sum_{i^k=1}^{m_{i^1...i^{k-1}}} a_{i^1...i^k} \le 1
\]
(when $k=1$ the end of the sum is $m$) and
\[
  0\le u_{i^1...i^k}<u_{i^1...i^k+1}<1\,, \quad u_{i^1...i^k+1}-u_{i^1...i^k}\ge a_{i^1...i^k}\,,
\]
when $k>1,\,i^k=m_{i^1...i^{k-1}}$ we substitute $u_{i^1...i^k+1}$ by 1. These hypotheses guarantee that the boxes
\[
     R_{i^1...i^d}=A_{i^1...i^d}([0,1]^d)
\]
have interiors that are pairwise disjoint, with edges parallel to the coordinate axes, the box $R_{i^1...i^d}$ having $k^{th}-$edge with length $a_{i^1...i^k}$. 
Geometrically, $\Lambda$ is constructed like the Sierpinski sponges, with the 1-approximation consisting of the boxes $R_{i^1...i^d}$, 
the 2-approximation consisting in replacing each box of the 1-aproximation by a rescaled affine copy of the 1-approximation, and so on. See Figure 1 for an illustration of the case $d=3$ (where we used $c_i,\,b_{ij},\,a_{ijk}$ instead of $a_{i^1},\,a_{i^1i^2},\,a_{i^1i^2i^3}$, respectively). When $d=2$ the sets $\Lambda$ are the Lalley-Gatzouras carpets, and their 1-approximation corresponds to the projection onto the $yz$-plane of Figure 1. In general we say that $\Lambda$ is a \emph{self-affine sponge}.\\\\

\begin{figure}[h]\label{fig}
\begin{center}
 \includegraphics[scale=0.45]{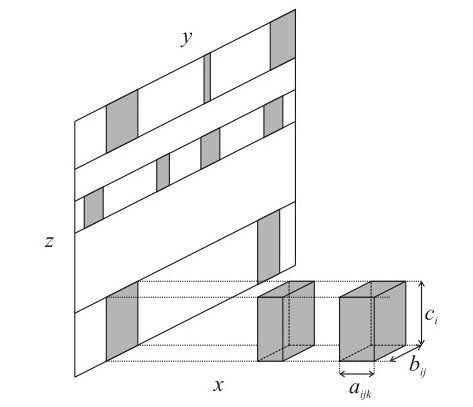}
  \caption{Self-affine sponges}
\end{center}
\end{figure}

Given ${\bf p} =(p_{i^1...i^{d-1}})$  a collection of
non-negative numbers satisfying
\begin{align*}
  &p_{i^1...i^k}=\sum_{i^{k+1}} p_{i^1...i^ki^{k+1}}\,,\quad k=d-2, d-3,...,1\\
  &\sum_{i^1} p_{i^1}=1,
\end{align*}
we define the numbers $\lambda_{k}({\bf{p}})$, $k=1,...,d-1$ by
\[
  \lambda_{k}({\bf{p}})=\frac{\underset{i^1...i^{k}}{\sum} p_{i^1...i^k}
  \log p_{i^1...i^k}-\underset{i^1...i^{k-1}}{\sum} p_{i^1...i^{k-1}}
 \log p_{i^1...i^{k-1}}}{\underset{i^1...i^{k}}{\sum} p_{i^1...i^k}\log
 a_{i^1...i^k}},
\]
(by convention: $0\log0=0$; for $k=1$ the second sum in the numerator is $0$)\\
and the number $t(\bf{p})$ as the unique real in $[0,1]$ satisfying
\[
  \sum_{i^1...i^{d-1}} p_{i^1...i^{d-1}} \log \left(\sum_{i^d}
  a_{i^1...i^{d-1}i^d}^{t(\bf{p})} \right)=0 .
\]
 The number $\sum_{k=1}^{d-1}\lambda_{k}({\bf{p}})$ is the Hausdorff dimension in the hyperplane $x_2...x_d$ of the set of generic points for the distribution ${\bf{p}}$; 
 the number $t(\bf{p})$ is the Hausdorff dimension of a typical $1$-dimensional fibre in the $x_1$-direction relative to the distribution ${\bf{p}}$, and is given by a 
 \emph{random Moran formula}. We will see that the sum of these two numbers is the Hausdorff dimension of a Bernoulli measure supported on $\Lambda$, and so we have the following result.

\begin{theorem}\label{lower}
Let $\Lambda$ be a self-affine sponge. Then
\begin{equation} \label{form}
 \hd \Lambda\ge\sup_{\bf{p}} \left\{ \sum_{k=1}^{d-1} \lambda_{k}({\bf{p}})+ t({\bf{p}}) \right\}.
\end{equation}
\end{theorem}
Here $\hd\Lambda$ stands for the Hausdorff dimension of a set $\Lambda$.

\vspace*{5pt}\noindent\textbf{Problem.} \emph{When does equality hold in (\ref{form})?} 
\vspace*{5pt}

In other words, when does the \emph{variational principle for the dimension} hold, in the class of self-affine sponges?
When $d=2$ this corresponds to Lalley-Gatzouras carpets and equality in (\ref{form}) was proved in \cite{10} (see also \cite{15} for an alternative proof and \cite{18} for a random version). 
For general $d$ but with $a_{i^1...i^{k}}=a^{(k)}$, $k=1,...,d$ this corresponds, essentially, to Sierpinski sponges and equality in (\ref{form}) was, essentially, proved in \cite{14}.

As we know now, by \cite{7},  equality in (\ref{form}) might not hold, in general, even when $d=3$.

For the sake of clarity we will restrict to $d=3$ (even though the results in this paper might extend to general $d$).  
In this case, we use the simpler notation $c_i$ for $a_i$, $i=1,...,m$ and $b_{ij}$ for $a_{ij}$, $j=1,...,m_j$, keeping the notation for $a_{ijk}$, $k=1,...,m_{ij}$ (see Figure 1). 
 
\section{$(a,b,c;\varepsilon)$-sponges}
In this section $d=3$.
\begin{definition} 
We say that a self-affine sponge is an $(a , b, c; \varepsilon)$-sponge, for some numbers $0 < a\le b \le c <1$, $\varepsilon\ge0$, if
\[
               e^{-\varepsilon}   \le  \frac{a_{ijk}}{a} \le e^\varepsilon, \quad e^{-\varepsilon}   \le  \frac{b_{ij}}{b} \le e^\varepsilon,
               \quad e^{-\varepsilon}   \le  \frac{c_{i}}{c} \le e^\varepsilon,
\]
for every $(i,j,k)\in \mathcal{I}$.             
\end{definition}

The case $\varepsilon=0$ corresponds essentially to what we called Sierpinski sponges (even though the numbers $a,b$ and $c$ need not be integers), and we will still call Sierpinski sponges to this larger class. As said before, their Hausdorff dimension was, essentially, computed in \cite{14} via the variational principle for the dimension. When we are \emph{close} to a Sierpinski sponge, say an $(a , b, c; \varepsilon)$-sponge with 
$\varepsilon>0$ small, we do not know if the variational principle for the dimension holds but we know that the Hausdorff dimension is close to the Hausdorff dimension of the Sierpinski sponge $(a , b, c; 0)$. When we talk about continuity of the Hausdorff dimension of self-affine sponges it is implicit that their alphabet $\mathcal{I}$ is fixed. 

Denote by $\Lambda_{a, b,c}$ a Sierpinski sponge $(a , b, c; 0)$, and by $\Lambda_{a, b,c; \varepsilon}$ an $(a, b,c; \varepsilon)-sponge$.
 
Let $\mathcal{J}=\{(i,j)\colon (i,j,k) \in\mathcal{I} \text{ for some } k\}$. Let
\[
\underline{t}=\min_{(i,j)\in\mathcal{J}} \frac{\log m_{ij}}{-\log a} \quad\text{and}\quad  
\overline{t}=\max_{(i,j)\in\mathcal{J}}\frac{\log m_{ij}}{-\log a}
\]

\begin{theorem}\label{cont}
Let $\Lambda_{a, b,c}$ be a Sierpinski sponge such that $\underline{t}<\overline{t}$. There exists a constant $C$ (depending only on $a$, $b$, $c$ and $\mathcal{I}$) such that for every $\varepsilon>0$ sufficiently small
\begin{equation*}
                                   \hd \Lambda_{a, b,c} -C \varepsilon     \le  \hd \Lambda_{a, b,c; \varepsilon}\le \hd \Lambda_{a, b,c} +C\varepsilon.
\end{equation*}                                   
\end{theorem}

\begin{rem}
Going through the proof of Theorem \ref{cont} it is possible to give an explicit expression for $C$, if that is needed for some application, but that is not our purpose in this paper.
\end{rem}

The lower estimate in Theorem \ref{cont} follows easily from Theorem \ref{lower}. For the upper estimate in Theorem \ref{cont} we will construct a 2-parameter family of Bernoulli measures.

Given $\Lambda$ a self-affine sponge we write
\[
   \mathrm{VP}(\Lambda)=\sup_{\bf{p}} \left\{ \lambda_{1}({\bf{p}})+\lambda_{2}({\bf{p}})+ t({\bf{p}}) \right\}.
\]   
Even though we do not know if the variational principle for the dimension holds for $(a,b,c;\varepsilon)$-sponges, $\varepsilon>0$ small, we have the following result.

\begin{corollary}\label{cor}
Let $\Lambda_{a, b,c}$ be a Sierpinski sponge such that $\underline{t}<\overline{t}$. There exists a constant $C$ (depending only on $a$, $b$, $c$ and $\mathcal{I}$) such that for every $\varepsilon>0$ sufficiently small
\begin{equation*}
                                 \mathrm{VP}(\Lambda_{a, b,c; \varepsilon}) \le  \hd \Lambda_{a, b,c; \varepsilon}\le \mathrm{VP}(\Lambda_{a, b,c; \varepsilon}) +C\varepsilon.
\end{equation*}                                   
\end{corollary}
The lower estimate in Corollary \ref{cor} follows from Theorem \ref{lower} and the upper estimate in Corollary \ref{cor} follows from the proof of Theorem \ref{cont}.

\section{Proof of Theorem \ref{lower}}

In this part, $d$ is any integer greater than or equal to 2.

There is a natural symbolic representation associated with our
system that we shall describe now. Consider the sequence space
$\Omega=\mathcal{I}^\mathbb{N}$. Elements of $\Omega$ will be
represented by $\omega=(\omega_1, \omega_2, ...)$ where
$\omega_n=(i^1_n,...,i^d_n)\in\mathcal{I}$. Given
$\omega\in\Omega$ and $n\in\mathbb{N}$, let $\omega(n)=(\omega_1,
\omega_2, ..., \omega_n)$ and define the \emph{cylinder of order
$n$},
\[
    C_{\omega(n)}=\{\omega'\in\Omega: \omega'_l=\omega_l,\,l=1,...,n\},
\]
and the \emph{basic box of order $n$},
\[
 R_{\omega(n)}=A_{\omega_1} \circ A_{\omega_2} \circ \cdots \circ A_{\omega_n}([0,1]^d).
\]
We have that $(R_{\omega(n)})_n$ is a decreasing sequence of
closed boxes having $k^{th}-$edge with length $\prod_{l=1}^n
a_{i_l^1...i_l^k}$. Thus $\bigcap_{n=1}^\infty R_{\omega(n)}$
consists of a single point which belongs to $\Lambda$ that we
denote by $\chi(\omega)$. This defines a continuous and surjective
map $\chi\colon \Omega\to\Lambda$ which is at most $2^d$ to 1, and
only fails to be a homeomorphism when some of the boxes
$R_{i^1...i^d}$ have nonempty intersection.

Let $\lambda({\bf{p}})=\sum_{k=1}^{d-1} \lambda_k({\bf{p}})$.
We shall construct probability measures $\mu_{\bf{p}}$ supported
on $\Lambda$ with
\[
\hd\mu_{\bf{p}}=\lambda({\bf{p}})+ t(\bf{p}).
\]
This gives what we want because $\hd\Lambda\ge\hd\mu_{\bf{p}}$.

Let $\tilde{\mu}_{\bf{p}}$ be the Bernoulli measure on $\Omega$
that assigns to each symbol $(i^1,...,i^d)\in\mathcal{I}$ the
probability
\[
    p_{i^1...i^{d-1}} \,\frac{a_{i^1...i^d}^{t(\bf{p})}}
    {\sum_{j^d} a_{i^1...i^{d-1}j^d}^{t(\bf{p})}}.
\]
In other words, we have
\[
    \tilde{\mu}_{\bf{p}}(C_{\omega(n)})=\prod_{l=1}^n p_{i_l^1...i_l^{d-1}}
    \,\frac{a_{i_l^1...i_l^d}^{t(\bf{p})}}{\sum_{j^d} a_{i_l^1...i_l^{d-1}j^d}^{t(\bf{p})}}.
\]
Let $\mu_{\bf{p}}$ be the probability measure on $\Lambda$ which
is the pushforward of $\tilde{\mu}_{\bf{p}}$ by $\chi$, i.e.
$\mu_{\bf{p}}=\tilde{\mu}_{\bf{p}}\circ\chi^{-1}$.

For calculating the Hausdorff dimension of $\mu_{\bf{p}}$, we
shall consider some special sets called \emph{approximate cubes}.
Given $\omega\in\Omega$ and $n\in\mathbb{N}$ such that $n\ge(\log
\min a_{i^1...i^d})/(\log \max a_{i^1})$, define
$L_n^0(\omega)=n$,
\begin{align}
 L_n^{1}(\omega)=&\max\left\{ k\ge1: \prod_{l=1}^n a_{i_l^1}\le
 \prod_{l=1}^k a_{i_l^1i_l^2}\right\}\notag\\
 &\vdots \label{ln2}\\
 L_n^{d-1}(\omega)=&\max\left\{ k\ge1: \prod_{l=1}^n a_{i_l^1}\le
 \prod_{l=1}^k a_{i_l^1...i_l^d}\right\}\notag
\end{align}
and the \emph{approximate cube}
\[
B_n(\omega)=\left\{
    \begin{matrix}
     \quad\quad \overline{i_{l}^1}=i_l^1,\: l=1,..., n\\
     \overline{\omega}\in\Omega:\:\, \overline{i_{l}^2}=i_l^2,\: l=1,...,L_n^1(\omega)\\
     \quad\quad\quad\quad \vdots\\
     \quad\quad\quad\quad\quad \overline{i_{l}^d}=i_l^d,\: l=1,...,L_n^{d-1}(\omega)
    \end{matrix}
   \right\}.
\]
We have that each approximate cube $B_n(\omega)$ is a finite union
of cylinder sets, and that approximate cubes are \emph{nested},
i.e., given two, say $B_n(\omega)$ and $B_{n'}(\omega')$, either
$B_n(\omega)\cap B_{n'}(\omega')=\emptyset$ or $B_n(\omega)\subset
B_{n'}(\omega')$ or $B_{n'}(\omega')\subset B_n(\omega)$.
Moreover, $\chi(B_n(\omega))=\tilde{B}_n(\omega)\cap \Lambda$
where $\tilde{B}_n(\omega)$ is a closed box in $\mathbb{R}^d$ with
edges parallel to the coordinate axes, the $k^{th}-$edge with
length $\prod_{l=1}^{L_n^{k-1}(\omega)} a_{i_l^1...i_l^k}$. By
(\ref{ln2}),
\begin{equation}\label{Q}
  1\le \frac{\underset{l=1}{\overset{L_n^{k-1}(\omega)}{\prod}} a_{i_l^1...i_l^k}}
  {\underset{l=1}{\overset{n}{\prod}} a_{i_l^1}}\le \max\,a_{i^1...i^d}^{-1},
\end{equation}
for $k=1,...,d$, hence the term ``approximate cube''. It follows
from (\ref{Q}) that
\begin{equation}\label{factor0}
  \frac{\sum_{l=1}^{L_n^{k-1}(\omega)} \log a_{i^1_l...i^k_l}}{\sum_{l=1}^n \log a_{i_l^1}}
  =1+\frac{1}{n} \frac{\sum_{l=1}^{L_n^{k-1}(\omega)} \log a_{i_l^1...i_l^k}-\sum_{l=1}^{n}
  \log a_{i_l^1}} {\frac{1}{n} \sum_{l=1}^{n}\log a_{i_l^1}}\to 1 .
\end{equation}
Also observe that $L_n^{k+1}(\omega)\le L_n^{k}(\omega)$ and
$L_n^k(\omega)\to\infty$ as $n\to\infty$.

First we calculate the dimension of the ``vertical'' part. Let
\[
 \mathcal{J}=\{(i^1,...,i^{d-1}): (i^1,...,i^{d-1},i^d)\in\mathcal{I} \text{ for some } i^d \}
\]
and $\Gamma=\mathcal{J}^\mathbb{N}$. Consider the natural
projections $\tilde{\pi}\colon\Omega\to\Gamma$ and
$\pi\colon\mathbb{R}^d\to\mathbb{R}^{d-1}$ given by
$\pi(x_1,...,x_d)=(x_2,...,x_d)$. We consider the measures
\[
  \tilde{\nu}_{\bf{p}}=\tilde{\mu}_{\bf{p}}\circ\tilde{\pi}^{-1}
  \quad \text{ and } \quad \nu_{\bf{p}}=\mu_{\bf{p}}\circ\pi^{-1}.
\]
\begin{lemma}\label{aux}
If $d>2$ then for every $k\in\{1,...,d-2\}$,
\[
  \frac{L_n^{k}(\omega)}{L_n^{k-1}(\omega)}\to\frac{\underset{i^1...i^{k}}{\sum} p_{i^1...i^{k}}
  \log a_{i^1...i^{k}}}{\underset{i^1...i^{k+1}}{\sum} p_{i^1...i^{k+1}} \log a_{i^1...i^{k+1}}}
  \, \text{ for $\tilde{\nu}_{\bf{p}}$-a.e. $\omega$}.
\]
\end{lemma}
\begin{proof}
It follows from (\ref{factor0}) that
\begin{equation}\label{factor1}
 \frac{\sum_{l=1}^{L_n^{k}(\omega)} \log a_{i^1_l...i^{k+1}_l}}
 {\sum_{l=1}^{L_n^{k-1}(\omega)} \log a_{i^1_l...i^{k}_l}}
 \to 1 .
\end{equation}
By Kolmogorov's Strong Law of Large Numbers (KSLLN),
\begin{equation}\label{factor2}
 \frac{1}{L_n^{k-1}(\omega)}\sum_{l=1}^{L_n^{k-1}(\omega)} \log
 a_{i^1_l...i^k_l}\to \underset{i^1...i^{k}}{\sum} p_{i^1...i^k}
 \log a_{i^1...i^k}\, \text{ for $\tilde{\nu}_{\bf{p}}$-a.e. $\omega$}
\end{equation}
and (redundantly)
\begin{equation}\label{factor3}
 \frac{1}{L_n^{k}(\omega)} \sum_{l=4}^{L_n^{k}(\omega)} \log
 a_{i^1_l...i^{k+1}_l}\to \underset{i^1...i^{k+1}}{\sum} p_{i^1...i^{k+1}}
 \log a_{i^1...i^{k+1}}\, \text{ for $\tilde{\nu}_{\bf{p}}$-a.e.
 $\omega$}.
\end{equation}
The result follows by (\ref{factor1}), (\ref{factor2}) and
(\ref{factor3}).
\end{proof}

The next lemma is a multidimensional version of \cite[Proposition 3.3]{10}.

\begin{lemma}\label{mult} $\hd\nu_{\bf{p}}=\lambda({\bf{p}})$. \end{lemma}
\begin{proof}
To calculate the Hausdorff dimension of $\nu_{\bf{p}}$ we are
going calculate its pointwise dimension and use \cite[Theorem 7.1]{21}.
Remember that $\chi(B_n(\omega))=\tilde{B}_n(\omega)\cap \Lambda$
where, by (\ref{Q}), $\tilde{B}_n(\omega)$ is ``approximately'' a
ball in $\mathbb{R}^d$ with radius $\prod_{l=1}^n a_{i_l^1}$, and
that
\[
   \nu_{\bf{p}}(\pi\tilde{B}_n(\omega))=\tilde{\nu}_{\bf{p}}(\tilde{\pi} B_n(\omega)).
\]
Also, $\chi$ is at most $2^d$ to 1. Taking this into account, by
 \cite[Theorem 7.1]{21} together with \cite[Theorem 15.3]{21}, one is left to prove that
\[
    \lim_{n\to\infty} \,\frac{\log \tilde{\nu}_{\bf{p}}(\tilde{\pi} B_n(\omega))}
    {\sum_{l=1}^n \log a_{i_l^1}}=\lambda({\bf{p}})
     \, \text{ for $\tilde{\nu}_{\bf{p}}$-a.e. $\omega$}.
\]
It follows from the definition of ${\tilde\nu}_{\bf{p}}$ that, for
${\tilde\nu}_{\bf{p}}$-a.e $\omega$, $p_{i_l^1...i_l^d}>0$ for
every $l$, so we may restrict our attention to these $\omega$. If
$d=2$ then $\tilde{\nu}_{\bf{p}}(\tilde{\pi}
B_n(\omega))=\prod_{l=1}^n p_{i_l^1}$ and the result follows by a
direct application of KSLLN. Otherwise we have that
\[
  \tilde{\nu}_{\bf{p}}(\tilde{\pi}B_n(\omega))=\prod_{l=1}^{L_n^{d-2}(\omega)}
  p_{i_l^1...i_l^{d-1}}\, \prod_{k=1}^{d-2}\,\prod_{l=L_n^{k}(\omega)+1}^{L_n^{k-1}(\omega)}
  p_{i_l^1...i_l^k}
\]
and
\begin{align*}
  \frac{\log \tilde{\nu}_{\bf{p}}(\tilde{\pi}B_n(\omega))}{\sum_{l=1}^n \log a_{i_l^1}}
  &= \frac{\sum_{l=1}^{L_n^{d-2}(\omega)} \log p_{i_l^1...i_l^{d-1}}}{\sum_{l=1}^{n}
   \log a_{i_l^1}}\\
   &\quad+ \sum_{k=1}^{d-2}\,\frac{\sum_{l=1}^{L_n^{k-1}(\omega)}
  \log a_{i_l^1...i_l^k}}{\sum_{l=1}^{n} \log a_{i_l^1}}\,
  \frac{\sum_{l=L_n^{k}(\omega)+1}^{L_n^{k-1}(\omega)} \log
  p_{i_l^1...i_l^k}} {\sum_{l=1}^{L_n^{k-1}(\omega)} \log
  a_{i_l^1...i_l^k}}\\
  &=\gamma_n+\sum_{k=1}^{d-2} \alpha_n^k\, \beta_n^k.
\end{align*}
By successive application of Lemma \ref{aux} one gets that
\[
    \frac{L_n^{d-2}(\omega)}{n}\to\frac{\underset{i^1}{\sum}\, p_{i^1} \log a_{i^1}}
    {\underset{i^1...i^{d-1}}{\sum} p_{i^1...i^{d-1}} \log a_{i^1...i^{d-1}}}
     \, \text{ for $\tilde{\nu}_{\bf{p}}$-a.e. $\omega$},
\]
and so, by KSLLN,
\[
   \gamma_n=\frac{\frac{L_n^{d-2}(\omega)}{n}\frac{1}{L_n^{d-2}(\omega)}
   \sum_{l=1}^{L_n^{d-2}(\omega)} \log p_{i_l^1...i_l^{d-1}}}{\frac{1}{n}\sum_{l=1}^{n}
   \log a_{i_l^1}}\to\frac{\underset{i^1...i^{d-1}}{\sum} p_{i^1...i^{d-1}}
  \log p_{i^1...i^{d-1}}}{\underset{i^1...i^{d-1}}{\sum} p_{i^1...i^{d-1}} \log a_{i^1...i^{d-1}}}
\]
for $\tilde{\nu}_{\bf{p}}$-a.e. $\omega$. By (\ref{factor0}),
$\alpha_n^k\underset{n\to\infty}{\to} 1$ for every $k$. We write
\[
  \beta_n^k = \frac{\frac{1}{L_n^{k-1}(\omega)} \sum_{l=1}^{L_n^{k-1}(\omega)} \log
  p_{i_l^1...i_l^k}-\frac{L_n^k(\omega)}{L_n^{k-1}(\omega)}
  \frac{1}{L_n^k(\omega)} \sum_{l=1}^{L_n^{k}(\omega)} \log
  p_{i_l^1...i_l^k}}{\frac{1}{L_n^{k-1}(\omega)} \sum_{l=1}^{L_n^{k-1}(\omega)} \log
  a_{i_l^1...i_l^k}}.
\]
Using KSLLN and Lemma \ref{aux} one gets that
\[
   \lim_{n\to\infty} \beta_n^k = \frac{\underset{i^1...i^{k}}{\sum} p_{i^1...i^k}\log
  p_{i^1...i^k}}{\underset{i^1...i^{k}}{\sum} p_{i^1...i^k}\log a_{i^1...i^k}}
  -\frac{\underset{i^1...i^{k}}{\sum} p_{i^1...i^k}\log
  p_{i^1...i^{k}}}{\underset{i^1...i^{k+1}}{\sum} p_{i^1...i^{k+1}}\log a_{i^1...i^{k+1}}}
\,\text{ for $\tilde{\nu}_{\bf{p}}$-a.e. $\omega$},
\]
and this gives what we want after a simple rearrangement.
\end{proof}

\begin{lemma}\label{dimmedida}
 $\hd \mu_{\bf{p}}=\lambda({\bf{p}})+ t(\bf{p})$.
\end{lemma}
\begin{proof}
As before, one is left to prove that
\[
    \lim_{n\to\infty} \,\frac{\log \tilde{\mu}_{\bf{p}}(B_n(\omega))}{\sum_{l=1}^n \log a_{i_l^1}}
    =\lambda({\bf{p}})+ t(\bf{p}) \, \text{ for $\tilde{\mu}_{\bf{p}}$-a.e. $\omega$}.
\]
We have that
\[
  \tilde{\mu}_{\bf{p}}(B_n(\omega))=\tilde{\nu}_{\bf{p}}(\tilde{\pi}B_n(\omega))
  \,\underbrace{\prod_{l=1}^{L_n^{d-1}(\omega)}\frac{a_{i_l^1...i_l^d}^{t(\bf{p})}}
  {\sum_{i^d} a_{i_l^1...i_l^{d-1}i^d}^{t(\bf{p})}}}_{\alpha_n}\,,
\]
By Lemma \ref{mult}, we only have to prove that
\[
    \lim_{n\to\infty} \,\frac{\log \alpha_n}{\sum_{l=1}^n \log a_{i_l^1}}
    =t(\bf{p}) \, \text{ for $\tilde{\mu}_{\bf{p}}$-a.e. $\omega$}.
\]
But
\begin{align*}
  \frac{\log \alpha_n}{\sum_{l=1}^n \log a_{i_l^1}}&=
   t({\bf{p}}) \, \frac{\sum_{l=1}^{L_n^{d-1}(\omega)} \log a_{i_l^1...i_l^d}}
   {\sum_{l=1}^{n} \log a_{i_l^1}}-\frac{\frac{1}{L_n^{d-1}(\omega)}
   \sum_{l=1}^{L_n^{d-1}(\omega)} \log \Bigr(\sum_{i^d} a_{i_l^1...i_l^{d-1}i^d}^{t(\bf{p})}
   \Bigl)}{\frac{n}{L_n^{d-1}(\omega)}\frac{1}{n}\sum_{l=1}^{n} \log a_{i_l^1}}\\
   &= t({\bf{p}})  \,  \beta_n \, - \,\frac{\overset{\,}{\gamma_n}}{\delta_n} \, .
\end{align*}
That $\beta_n\to 1$ follows from (\ref{factor0}). Now we can write
\[
  \gamma_n = \sum_{i^1...i^{d-1}} \frac{P(\omega,L_n^{d-1}(\omega),i^1...i^{d-1})}
  {L_n^{d-1}(\omega)}\,\log \left(\sum_{i^d} a_{i^1...i^{d-1}i^d}^{t(\bf{p})} \right),
\]
where
\[
  P(\omega,n,i^1...i^{d-1})=\sharp \{1\le l\le n: (i_l^1...i_l^{d-1})=(i^1...i^{d-1})\}\,
\]
for $(i^1,...,i^{d-1})\in\mathcal{J}$. By KSLLN,
\[
   \frac{P(\omega,n,i^1...i^{d-1})}{n}\to p_{i^1...i^{d-1}} \,
   \text{ for $\tilde{\mu}_{\bf{p}}$-a.e. $\omega$},
\]
so, by the definition of $t(\bf{p})$,
\[
   \gamma_n\to 0\,\text{ for $\tilde{\mu}_{\bf{p}}$-a.e. $\omega$}.
\]
Since $n/L_n^{d-1}(\omega)\ge 1$, we have that $| \delta | \ge
\log \,(\min\, a_{i^1}^{-1})>0$, so we also have that
\[
   \frac{\gamma_n}{\delta_n} \to 0\,\text{ for $\tilde{\mu}_{\bf{p}}$-a.e. $\omega$},
\]
thus completing the proof.
\end{proof}
As noticed in the beginning of this section, these lemmas imply
\[
   \hd\Lambda\ge\sup_{\bf{p}} \left\{ \lambda({\bf{p}})+t(\bf{p})\right\}.
\]

\begin{rem}
Theorem \ref{lower} can be extended to \emph{non-linear self-affine sponges} by using the bounded distortion techniques employed in \cite{15} (see also \cite{19}).
\end{rem}

\section{A 2-parameter family of Bernoulli measures}
In this section $d=3$ and $\Lambda$ is a self-affine sponge.

We start by generalizing the definitions of $\underline{t}$ and $\overline{t}$ for any self-affine sponge.
Let $\underline{t}=\min_{{\bf{p}}} t({\bf{p}})$ and $\overline{t}=\max_{{\bf{p}}} t({\bf{p}})$. 
For $(i,j)\in\mathcal{J}$, define $t_{ij}$ to be the unique real in $[0,1]$ satisfying
\begin{equation}\label{tij}
    \sum_{k=1}^{m_{ij}} a_{ijk}^{t_{ij}}=1.
\end{equation}
It is easy to see that 
\[
\underline{t}=\min_{(i,j)\in\mathcal{J}} t_{ij} \quad\text{and}\quad  \overline{t}=\max_{(i,j)\in\mathcal{J}} t_{ij}
\]
We observe that the condition $\underline{t}<\overline{t}$ is \emph{open} in the numbers $a_{ijk}$, so if it is satisfied for an
$(a,b,c;0)$-sponge then it is also satisfied for an $(a,b,c;\varepsilon)$-sponge, for some $\varepsilon=\varepsilon(a, \mathcal{I})>0$.

Let
\[
 \mathcal{P}=\left\{ {\bf{p}}=(p_{ij})_{(i,j)\in\mathcal{J}} : p_{ij}>0 \text{ for all }
 (i,j)\in\mathcal{J} \text{ and } \sum_{i,j} p_{ij}=1\right\}.
\]

We will need the following \emph{generic} hypothesis on the numbers $a_{ijk}$. For each
$t\in (\underline{t},\overline{t})$, there exist $1\le i\le m$ and $1\le j<j'\le m_i$
such that
\begin{equation}\label{hip}
  \sum_{k=1}^{m_{ij}} a_{ijk}^t \ne \sum_{k=1}^{m_{ij'}} a_{ij'k}^t.
\end{equation}

In the next lemma there will be no restrictions on the numbers $a_{ijk}$ (beside hypothesis (\ref{hip}). In this way, we say that a self-affine sponge is an $(c,b;\varepsilon)$-sponge, for some numbers $0 < b \le c <1$, $\varepsilon\ge0$, if
\[
               e^{-\varepsilon}   \le  \frac{b_{ij}}{b} \le e^\varepsilon,
       \quad e^{-\varepsilon}   \le  \frac{c_{i}}{c} \le e^\varepsilon,
\]
for every $(i,j)\in \mathcal{J}$.             

\begin{lemma}\label{lemapij}
Let $0<b\le c\le1$, $0<a_{ij}\le1$, $(i,j)\in\mathcal{J}$ and $\delta>0$. There exists $\varepsilon=\varepsilon(c,b, a_{ij},\mathcal{I},\delta)$ such that if $\Lambda$ is an $(b,c; \varepsilon)$-sponge satisfying (\ref{hip}) and $\underline{t}<\overline{t}$ then:  given $t\in [\underline{t}+\delta, \overline{t}-\delta]$ and $\rho\in [\delta,1]$,
there exists a probability vector ${\bf{p}}={\bf{p}}(t,\rho)$ such that $t({\bf{p}})=t$ and
\begin{equation*}
  {p}_{ij}= c_i^{\lambda_1({\bf{p}})}\, b_{ij}^{\lambda_2({\bf{p}})} \Bigl(\sum_{k} a_{ijk}^{t}\Bigr)^{\alpha}\,
  \Bigl(\sum_j b_{ij}^{\lambda_2({\bf{p}})} \Bigl(\sum_{k} a_{ijk}^{t}\Bigr)^{\alpha}\Bigr)^{\rho-1}, \quad (i,j)\in\mathcal{J},
\end{equation*}
where $\alpha=\alpha(t,\rho)\in\mathbb{R}$ is a $\mathrm{C}^1$ function defined in $[\underline{t}+\delta,\overline{t}-\delta]\times[\delta,1]$.  
Moreover, for each $\rho\in (0,1]$,
$\alpha(\underline{t}+\delta,\rho)\to-\infty$ and $\alpha(\overline{t}-\delta,\rho)\to\infty$ when $\delta\to0$.
\end{lemma}
\begin{proof}
Given $\alpha,\lambda_1\in\mathbb{R}$, $t\in
(\underline{t}, \overline{t})$, $\rho\in (0,1]$ and $\lambda_2\in[0,1]$, we define a
probability vector ${\bf{p}}(\alpha,\lambda_1,\lambda_2,t,\rho)$
by
\begin{equation*}\label{pij0}
  p_{ij}(\alpha,\lambda_1,\lambda_2,t,\rho)= C(\alpha,\lambda_1,\lambda_2,t,\rho)\,
  c_i^{\lambda_1}\, b_{ij}^{\lambda_2} \Bigl(\sum_{k} a_{ijk}^{t}\Bigr)^{\alpha}\,
  \gamma_i(\alpha,\lambda_2,t)^{\rho-1}
\end{equation*}
where
\[
  \gamma_i(\alpha,\lambda_2,t)= \sum_j b_{ij}^{\lambda_2}\Bigl(\sum_{k} a_{ijk}^{t}\Bigr)^{\alpha}
\]
and
\[
  C(\alpha,\lambda_1,\lambda_2,t,\rho)=\Bigl(\sum_i c_i^{\lambda_1}\,
  \gamma_i(\alpha,\lambda_2,t)^\rho\Bigl)^{-1},
\]
for each $(i,j)\in\mathcal{J}$.

Let $F$ be the continuous function defined by
\[
   F(\alpha,\lambda_1,\lambda_2,t,\rho)=\sum_{i,j} p_{ij}(\alpha,\lambda_1,\lambda_2,t,\rho)
   \log\Bigl(\sum_k a_{ijk}^t\Bigr).
\]
We are going to prove there exists a unique
$\alpha=\alpha(\lambda_1,\lambda_2,t,\rho)$, continuously varying,
such that $F(\alpha,\lambda_1,\lambda_2,t,\rho)=0$, i.e.
$t({\bf{p}}(\alpha,\lambda_1,\lambda_2,t,\rho))=t$.

\emph{Unicity.} We have that, for each $(i,j)\in\mathcal{J}$,
\[
  \frac{\partial p_{ij}}{\partial \alpha}=\frac{1}{C} \frac{\partial C}{\partial \alpha} p_{ij}
  +\log\Bigl(\sum_k a_{ijk}^t\Bigr) p_{ij}+(\rho-1)\frac{1}{\gamma_i}
  \frac{\partial \gamma_{i}}{\partial \alpha} p_{ij}.
\]
Also,
\begin{equation}\label{derivative1}
  \frac{1}{C}\frac{\partial C}{\partial \alpha}=-\rho\sum_i p_i \frac{1}{\gamma_i}
  \frac{\partial \gamma_{i}}{\partial \alpha}
\end{equation}
and
\begin{equation}\label{derivative2}
 \frac{1}{\gamma_i} \frac{\partial \gamma_{i}}{\partial \alpha}=
 \sum_j \frac{p_{ij}}{p_i} \log\Bigl(\sum_k a_{ijk}^t\Bigr)
\end{equation}
where
\begin{equation}\label{pi}
  p_i=\sum_j p_{ij}=C c_i^{\lambda_1} \gamma_i^\rho.
\end{equation}
So, by simple rearrangement we get
\begin{align*}\label{Falfa}
 \frac{\partial F}{\partial \alpha}&=\sum_{i,j} \frac{\partial p_{ij}}{\partial \alpha}
   \log\Bigl(\sum_k a_{ijk}^t\Bigr)\\ 
   &=\rho\left\{\sum_i p_i \Bigl(\sum_j \frac{p_{ij}}{p_i}
   \log\Bigl(\sum_k a_{ijk}^t\Bigr)\Bigr)^2-\Bigl(\sum_{i,j}p_{ij}
   \log\Bigl(\sum_k a_{ijk}^t\Bigr)\Bigr)^2\right\}  \notag \\
   &\quad+\sum_i p_i \left\{ \sum_j \frac{p_{ij}}{p_i} \Bigl(\log\Bigl(\sum_k a_{ijk}^t\Bigr)\Bigr)^2
   -\Bigl(\sum_j \frac{p_{ij}}{p_i} \log\Bigl(\sum_k a_{ijk}^t\Bigr)\Bigr)^2\right\}.\notag
\end{align*}
By the Cauchy-Schwarz inequality we have that the expressions
between curly brackets are non-negative and the second one is
positive if there exists $i\in\{1,...,m\}$ such that the function
\[
   j\mapsto \sum_k a_{ijk}^t
\]
is non-constant (note that ${\bf{p}}\in\mathcal{P}$). This is guaranteed by hypothesis
(\ref{hip}). Thus $\partial F/\partial \alpha>0$.

\emph{Existence.} For fixed $(\lambda_1,\lambda_2,t,\rho)$, we
will look at the limit distributions of
${\bf{p}}(\alpha)={\bf{p}}(\alpha,\lambda_1,\lambda_2,t,\rho)$ as
$\alpha$ goes to $+\infty$ and $-\infty$. 
Let
\[
    A_t=\max_{(i,j)\in\mathcal{J}} \sum_k a_{ijk}^t.
\]
For $(i,j)\in\mathcal{J}$ such that $t<t_{ij}$ (remember the definition (\ref{tij})) we have that
\[
   \sum_k a_{ijk}^t>\sum_k a_{ijk}^{t_{ij}}=1,
\]
so $A_t>1$. Consider $(\bar{i},\bar{j}),(i,j)\in\mathcal{J}$ such
that
\[
   \sum_k a_{ijk}^t<\sum_k a_{\bar{i}\bar{j}k}^t=A_t.
\]
We have
\[
    \frac{p_{ij}(\alpha)}{p_{\bar{i}\bar{j}}(\alpha)}\le D
    \left(\frac{\sum_k a_{ijk}^t}{\sum_k a_{\bar{i}\bar{j}k}^t}\right)^\alpha
    \left(\frac{\gamma_{\bar{i}}(\alpha)}{\gamma_{i}(\alpha)}\right)^{1-\rho},
\]
for some constant $D$ not depending on $\alpha$. Now, for
$\alpha>0$,
\[
 \frac{\gamma_{\bar{i}}(\alpha)}{\gamma_{i}(\alpha)}\le \tilde{D}
 \left(\frac{\sum_k a_{\bar{i}\bar{j}k}^t}{\sum_k a_{ijk}^t}\right)^{\alpha},
\]
for some constant $\tilde{D}$ not depending on $\alpha$. So,
\begin{equation*}
 \frac{p_{ij}(\alpha)}{p_{\bar{i}\bar{j}}(\alpha)}\le D \tilde{D}
  \left(\frac{\sum_k a_{ijk}^t}{\sum_k a_{\bar{i}\bar{j}k}^t}\right)^{\alpha\rho},
\end{equation*}
which converges to 0 as $\alpha\to\infty$. This implies that
\begin{equation}\label{exist1}
  \sum_{i,j} p_{ij}(\alpha)\log\Bigl(\sum_k a_{ijk}^t\Bigr)
  \underset{\alpha\to\infty}{\longrightarrow}\log A_t>0.
\end{equation}
In the same way, defining
\[
    B_t=\min_{(i,j)\in\mathcal{J}} \sum_k a_{ijk}^t<1,
\]
and taking $(\underline{i},\underline{j}),(i,j)\in\mathcal{J}$
such that
\[
   B_t=\sum_k a_{\underline{i}\underline{j}k}^t<\sum_k a_{ijk}^t,
\]
we get, for $\alpha<0$,
\begin{equation*}
 \frac{p_{ij}(\alpha)}{p_{\underline{i}\underline{j}}(\alpha)}\le D \tilde{D}
  \left(\frac{\sum_k a_{ijk}^t}{\sum_k a_{\underline{i}\underline{j}k}^t}\right)^{\alpha\rho},
\end{equation*}
which converges to 0 as $\alpha\to-\infty$. This implies that
\begin{equation}\label{exist2}
  \sum_{i,j} p_{ij}(\alpha)\log\Bigl(\sum_k a_{ijk}^t\Bigr)
  \underset{\alpha\to-\infty}{\longrightarrow}\log B_t<0.
\end{equation}
By (\ref{exist1}), (\ref{exist2}) and continuity, there exists
$\alpha\in\mathbb{R}$ such that
$F(\alpha,\lambda_1,\lambda_2,t,\rho)=0$. The continuity of
$\alpha(\lambda_1,\lambda_2,t,\rho)$ follows from the uniqueness
part and the implicit function theorem. Actually, since
$F(\alpha,\lambda_1,\lambda_2,t,\rho)$ is continuously differentiable, we also
get that $\alpha(\lambda_1,\lambda_2,t,\rho)$ is continuoulsly differentiable. 
Observe that
$t({\bf{p}})=\overline{t}\Rightarrow {\bf{p}}\in\partial\mathcal{P}$ (in this lemma we are
assuming $\underline{t}<\overline{t}$), so since
\[
    t({\bf{p}}(\alpha(\lambda_1,\lambda_2,t,\rho)))\to \overline{t} \quad\text{when}\quad
    t\to\overline{t}
\]
then
\[
   {\bf{p}}(\alpha(\lambda_1,\lambda_2,t,\rho))\to\partial\mathcal{P}
   \quad\text{when}\quad t\to\overline{t},
\]
which implies
\[
  \alpha(\lambda_1,\lambda_2,t,\rho)\to\infty \quad\text{when}\quad t\to\overline{t}
\]
(this convergence is uniform in $\lambda_1,\lambda_2\in [0,1]$, and $\rho$ in a compact set set of $(0,1]$).
In the same way we see that
\[
  \alpha(\lambda_1,\lambda_2,t,\rho)\to-\infty \quad\text{when}\quad t\to\underline{t}.
\]
Moreover,
\begin{equation}\label{alfalambda1}
 \frac{\partial \alpha}{\partial \lambda_1}=-\left(\frac{\partial F}{\partial \alpha}\right)^{-1} \frac{\partial F}{\partial \lambda_1},
\end{equation}
where
\begin{align}\label{Flambda1}
\frac{\partial F}{\partial \lambda_1}=\sum_{i,j} \frac{\partial p_{ij}}{\partial \lambda_1} \log \Bigl( \sum_k a_{ijk}^t \Bigr)
   =\sum_{i,j} p_{ij} \log c_i \log \Bigl( \sum_k a_{ijk}^t \Bigr),
\end{align}
and
\begin{equation}\label{alfalambda2}
 \frac{\partial \alpha}{\partial \lambda_2}=-\left(\frac{\partial F}{\partial \alpha}\right)^{-1} \frac{\partial F}{\partial \lambda_2},
\end{equation}
where
\begin{align}\label{Flambda2}
\frac{\partial F}{\partial \lambda_2}=\sum_{i,j} \frac{\partial p_{ij}}{\partial \lambda_2} \log \Bigl( \sum_k a_{ijk}^t \Bigr)
   =\sum_{i,j} p_{ij}\left(  \log b_{ij} + (\rho-1) \frac{1}{\gamma_i} \frac{\partial \gamma_i}{\partial \lambda_2} \right) \log \Big( \sum_k a_{ijk}^t \Bigr),
\end{align}
and
\begin{equation}\label{gammalambda2}
 \frac{1}{\gamma_i} \frac{\partial \gamma_i}{\partial \lambda_2}= \sum_{j} \frac{p_{ij}}{p_i}  \log b_{ij}.
\end{equation}

Now we want to find $\lambda_1=\lambda_1(\lambda_2,t,\rho)$
continuoulsly differentiable such that
\begin{equation}\label{fixc}
   C(\alpha(\lambda_1,\lambda_2,t,\rho),\lambda_1,\lambda_2,t,\rho)=1.
\end{equation}
We have
\[
  \frac{\partial}{\partial \lambda_1}
  C(\alpha(\lambda_1,\lambda_2,t,\rho),\lambda_1,\lambda_2,t,\rho)=
  \frac{\partial C}{\partial \alpha} \frac{\partial \alpha}{\partial \lambda_1}+
  \frac{\partial C}{\partial \lambda_1}
\]
Observe that, by (\ref{derivative1}) and (\ref{derivative2}),
\begin{equation*}
   \frac{\partial C}{\partial \alpha}=0
\end{equation*}
(at points
$(\alpha(\lambda_1,\lambda_2,t,\rho),\lambda_1,\lambda_2,t,\rho)$),
and
\begin{equation}\label{derivC1}
  \frac{\partial \log C}{\partial \lambda_1}=\frac{1}{C}\frac{\partial C}{\partial \lambda_1}
  =- \sum_i p_i \log c_i \ge \min_i \log c_i^{-1}>0.
\end{equation}
So,
$C(\alpha(\lambda_1,\lambda_2,t,\rho),\lambda_1,\lambda_2,t,\rho)$
is a continuoulsly differentiable function, for each $(\lambda_2,t,\rho)$, is
strictly increasing in $\lambda_1$ and, by (\ref{derivC1}), has
limit $\infty$ as $\lambda_1\to\infty$ and limit 0 as
$\lambda_1\to-\infty$. By the implicit function theorem, there is
a unique $\lambda_1=\lambda_1(\lambda_2,t,\rho)$, which is continuoulsly
differentiable, satisfying (\ref{fixc}). 
Moreover,
\begin{equation}\label{derivC2}
 \frac{\partial \lambda_1}{\partial \lambda_2}=
  -\left( \frac{\partial C}{\partial \lambda_1}\right)^{-1} \frac{\partial C}{\partial \lambda_2}
\end{equation}
and
\begin{equation}\label{derivC3}
  \frac{1}{C} \frac{\partial C}{\partial \lambda_2}
   =-\rho \sum_i p_i \frac{1}{\gamma_i}\frac{\partial \gamma_i}{\partial \lambda_2}.
\end{equation}
So, by (\ref{derivC1})-(\ref{derivC3}) and (\ref{gammalambda2}), we get
\begin{equation}\label{derivC5}
 \frac{\partial \lambda_1}{\partial \lambda_2}=
 -\rho\frac{\sum_{i,j} p_{ij}\log b_{ij}}{\sum_i p_i \log c_i}
\end{equation}

We see that
\begin{align*}
 \lambda_1({\bf{p}})&=\lambda_1+\rho \,
 \frac{\sum_i p_i \log \gamma_i}
 {\sum_i p_i \log c_i},\notag \\ 
 \lambda_2({\bf{p}})&=\lambda_2-
 \frac{\sum_i p_i \log \gamma_i}
 {\sum_{i,j} p_{ij} \log b_{ij}} \label{gama}.
\end{align*}

We use the following notation
\begin{align*}
   \tilde\Theta(\lambda_2,t,\rho)&=(\alpha(\lambda_1(\lambda_2,t,\rho),\lambda_2,t,\rho),
   \lambda_1(\lambda_2,t,\rho),\lambda_2,t,\rho),\\
   \Theta(\lambda_2,t,\rho)&=(\alpha(\lambda_1(\lambda_2,t,\rho),\lambda_2,t,\rho),
   \lambda_2,t,\rho).
\end{align*}
We want to prove there exists a unique $\lambda_2=\lambda_2(t,\rho)$, continuously differentiable, such that
\begin{equation*}\label{gama0}
   H(\lambda_2,t,\rho) := \sum_i p_i(\tilde\Theta)\log \gamma_i(\Theta)=0.
\end{equation*}
By (\ref{pi}) and (\ref{fixc}),
\begin{align*}
  &\frac{\partial}{\partial \lambda_2} \sum_i p_i(\tilde\Theta)\log \gamma_i(\Theta)=\frac{\partial \lambda_1}{\partial \lambda_2} \sum_i      p_i(\tilde\Theta)\log c_i \log \gamma_i(\Theta)\\
  &+\rho\sum_i p_i(\tilde\Theta)\frac{1}{\gamma_i(\Theta)}\frac{\partial}{\partial \lambda_2} 
  \gamma_i(\Theta) \log \gamma_i(\Theta)+ \sum_i p_i(\tilde\Theta)\frac{1}{\gamma_i(\Theta)}\frac{\partial}{\partial \lambda_2} 
  \gamma_i(\Theta)\notag\\
  \end{align*}
We have that
\begin{equation*}\label{gama2}
  \frac{1}{\gamma_i(\Theta)}\frac{\partial}{\partial \lambda_2}\gamma_i(\Theta)=
  \frac{1}{\gamma_i(\Theta)}\frac{\partial\gamma_i}{\partial \alpha}(\Theta)
  \left(\frac{\partial \alpha}{\partial\lambda_1}\frac{\partial\lambda_1}{\partial \lambda_2}
  +\frac{\partial\alpha}{\partial\lambda_2}\right)+
  \frac{1}{\gamma_i(\Theta)}\frac{\partial\gamma_i}{\partial \lambda_2}(\Theta).
\end{equation*}
Then, using (\ref{derivC5}) and (\ref{gammalambda2}), we get 
\begin{align}
  &\frac{\partial}{\partial \lambda_2} \sum_i p_i(\tilde\Theta)\log \gamma_i(\Theta)=\sum_{i,j} p_{ij}(\tilde\Theta)\log b_{ij}\notag\\
  &+\rho \sum_i p_i(\tilde\Theta) \left\{  \sum_j \frac{p_{ij}(\tilde\Theta)}{p_{i}(\tilde\Theta)} \log b_{ij}-\log c_i \,\frac{\sum_{i,j} p_{ij}(\tilde\Theta)\log b_{ij}}{\sum_{i} p_{i}(\tilde\Theta)\log c_{i}} \right\} \log \gamma_i(\Theta)\label{curly}\\
  &+\rho \left(\frac{\partial \alpha}{\partial\lambda_1}\frac{\partial\lambda_1}{\partial \lambda_2}
  +\frac{\partial\alpha}{\partial\lambda_2}\right) \sum_i p_i(\tilde\Theta) \frac{1}{\gamma_i(\Theta)}\frac{\partial\gamma_i}{\partial \alpha}(\Theta) \log \gamma_i(\Theta). \label{fodao}
  \end{align}
The term (\ref{curly}) can be made arbitrarily small if $c_i$ and $b_{ij}$ are sufficiently close to $c$ and $b$, respectively. Note that $\log \gamma_i(\Theta)$ is uniformly bounded for $t\in[\underline{t}+\delta, \overline{t}-\delta]$ and $\rho\in[\delta,1]$. The term (\ref{fodao}) can also be made arbitrarily small if $c_i$ and $b_{ij}$ are sufficiently close to $c$ and $b$, respectively, because, by (\ref{alfalambda1})-(\ref{gammalambda2}), $\frac{\partial \alpha}{\partial\lambda_1}$ and $\frac{\partial \alpha}{\partial\lambda_2}$ satisfy this property (and the other quantities are uniformly bounded for $t\in[\underline{t}+\delta, \overline{t}-\delta]$ and $\rho\in[\delta,1]$, see (\ref{derivative2})).
Then
\[
\frac{\partial H}{\partial \lambda_2}=\frac{\partial}{\partial \lambda_2} \sum_i p_i(\tilde\Theta)\log \gamma_i(\Theta)<0
\]
if $\Lambda$ is a $(c,b;\varepsilon)$-sponge, for some  $\varepsilon=\varepsilon(c,b,a_{ij},\mathcal{I},\delta)>0$. Then the existence of $\lambda_2(t,\rho)$ as claimed follows from the implicit function theorem. 
\end{proof}

\section{Proof of Theorem \ref{cont} and Corollary \ref{cor}}

It follows from definitions (\ref{ln2}) that, for $(a,b,c; \varepsilon)$-sponges,
\begin{equation*}
                           -\frac{1}{n}    +   \frac{L_n^1 \varepsilon}{n \log b}    + \frac{\varepsilon}{\log b} +     \frac{\log c}{\log b}   \le
  \frac{L_n^1}{n}\le \frac{\log c}{\log b} - \frac{\varepsilon}{\log b} - \frac{L_n^1 \varepsilon}{n \log b},
\end{equation*}
and similarly,  
\begin{equation*}
   -\frac{1}{n}    +   \frac{L_n^2 \varepsilon}{n \log a}    + \frac{\varepsilon}{\log a} +     \frac{\log c}{\log a}   \le
  \frac{L_n^2}{n}\le \frac{\log c}{\log a} - \frac{\varepsilon}{\log a} - \frac{L_n^2 \varepsilon}{n \log a}.
\end{equation*}
And so,
\begin{equation}\label{est1}
                           -\frac{1}{n} -  A \varepsilon +     \frac{\log c}{\log b}   \le \frac{L_n^1}{n}\le \frac{\log c}{\log b} + A\varepsilon,
\end{equation}
and 
\begin{equation}\label{est2}
                           -\frac{1}{n} -  A \varepsilon +     \frac{\log c}{\log a}   \le \frac{L_n^2}{n}\le \frac{\log c}{\log a} + A\varepsilon,
\end{equation}
where $A$ is a positive constant (depending only on $a,b,c$) and for every $\varepsilon>0$ sufficiently small.

We begin by proving the lower estimate in Theorem \ref{cont}.
We leave to the reader to prove that, for $(a,b,c; \varepsilon)$-sponges,
\begin{align*}
           &\frac{\sum_{i=1}^m p_{i} \log p_{i}}{\sum_{i=1}^m p_{i}\log c_{i}} + 
            \frac{\sum_{i=1}^m \sum_{j=1}^{m_i} p_{ii} \log p_{ij} - \sum_{i=1}^m p_{i} \log p_{i}}{\sum_{i=1}^m \sum_{j=1}^{m_i} p_{ij}\log b_{ij}} \\
            &\ge \frac{\sum_{i=1}^m p_{i} \log p_{i}}{\log c} + 
            \frac{\sum_{i=1}^m \sum_{j=1}^{m_i} p_{ii} \log p_{ij} - \sum_{i=1}^m p_{i} \log p_{i}}{\log b}-B\varepsilon,
\end{align*}
where $B$ is a positive constant (depending only on $a,b,c$ and $\mathcal{I}$) and for every $\varepsilon>0$ sufficiently small.
Let $t_0$ be such that
\[
    \sum_{i=1}^m \sum_{j=1}^{m_i} p_{ij} \log \Bigl( \sum_{k=1}^{m_{ij}} a^{t_0} \Bigr)=0.
\]
We want to see that $t({\bf{p}})\ge t_0 -D\varepsilon$ for some positive constant $D$ (depending only on $a,b,c$ and $\mathcal{I}$) and for every $\varepsilon>0$ sufficiently small. This is true if
\[
    \sum_{i=1}^m \sum_{j=1}^{m_i} p_{ij} \log \Bigl( \sum_{k=1}^{m_{ij}} a_{ijk}^{t_0-D\varepsilon} \Bigr)\ge 0.
\]
Now
\begin{align*}
    &\sum_{i=1}^m \sum_{j=1}^{m_i} p_{ij} \log \Bigl( \sum_{k=1}^{m_{ij}} a_{ijk}^{t_0-D\varepsilon} \Bigr) \ge
    \sum_{i=1}^m \sum_{j=1}^{m_i} p_{ij} \log \Bigl( \sum_{k=1}^{m_{ij}} (ae^{-\varepsilon})^{t_0-D\varepsilon} \Bigr)\\
    &\ge D\varepsilon \log a^{-1}-\varepsilon (t_0-D\varepsilon)\ge 0,
\end{align*}
for some appropriate $D$ as before and for every $\varepsilon>0$ sufficiently small. Then we have that
$\mathrm{VP}(\Lambda_{a,b,c;\varepsilon})\ge\mathrm{VP}(\Lambda_{a,b,c})-C\varepsilon$, for some positive constant $C$ (depending only on $a,b,c$ and $\mathcal{I}$) and for every $\varepsilon>0$ sufficiently small.
Using Theorem \ref{lower} and that $\hd \Lambda_{a,b,c} =\mathrm{VP} (\Lambda_{a,b,c})$ we get
$\hd \Lambda_{a,b,c;\varepsilon} \ge \hd \Lambda_{a,b,c}-C\varepsilon$, as we wish.

\begin{rem}
The fact that $\hd \Lambda_{a,b,c} =\mathrm{VP} (\Lambda_{a,b,c})$ was \emph{essentially} proved in \cite{11}. Although
in \cite{11} the numbers $a^{-1}, b^{-1}, c^{-1}$ are assumed to be integers, the proofs work the same when these numbers are not integers. Alternatively, this paper (see the next two lemmas) with $\varepsilon=0$ also gives this result.
\end{rem}

Now we prove the upper estimate in Theorem \ref{cont}.
\begin{lemma}\label{lemasup}
Let $0<a\le b\le c\le1$ and assume $\underline{t}<\overline{t}$. There exists a positive constant $\tilde{C}$  (depending only on $a,b,c$ and $\mathcal{I}$) such that for every $\varepsilon>0$ sufficiently small we have the following: if $\Lambda$ is a $(a,b,c; \varepsilon)$-sponge and $\omega\in\Omega$ there exists ${\bf{p}}\in\mathcal{P}$ such that
\[
 \liminf_{n\to\infty} \frac{\log \tilde{\mu}_{\bf{p}}(B_n(\omega))}
 {\sum_{l=1}^n \log c_{i_l}}\le \mathrm{VP}(\Lambda)+\tilde{C}\varepsilon.
\]
\end{lemma}
\begin{proof}
First assume that $\Lambda$ satisfy hypothesis (\ref{hip}) (at the end of the proof we say how to deal with the general case). Fix $\omega\in\Omega$. We use the notation
\[
   d_{{\bf{p}},n}(\omega)=\frac{\log \tilde{\mu}_{{\bf{p}}}(B_n(\omega))}
   {\sum_{l=1}^n \log c_{i_l}}.
\]
Then it follows from the proofs of Lemma \ref{mult} and Lemma \ref{dimmedida} that, if
${\bf{p}}\in\mathcal{P}$,
\begin{align}
  d_{{\bf{p}},n}(\omega)=&
  \frac{\sum_{l=1}^n \log p_{i_l}}{\sum_{l=1}^n \log c_{i_l}}+\beta_n(\omega)\,
  \frac{\sum_{l=1}^{L_n^1(\omega)}\log p_{i_lj_l}-\sum_{l=1}^{L_n^1(\omega)}
  \log p_{i_l}}{\sum_{l=1}^{L_n^1(\omega)} \log b_{i_lj_l}} \label{vert1}\\
  &+\eta_n(\omega) t({\bf{p}})-\frac{\sum_{l=1}^{L_n^2(\omega)}
  \log\Bigl(\sum_k a_{i_lj_lk}^{t({\bf{p}})}\Bigr)}{\sum_{l=1}^n \log c_{i_l}}\notag
\end{align}
where, by (\ref{factor0}),
\[
  \beta_n(\omega)=\frac{\sum_{l=1}^{L_n^1(\omega)} \log b_{i_lj_l}}{\sum_{l=1}^n \log c_{i_l}}
  \underset{n\to\infty}{\longrightarrow} 1
\]
and
\[
  \eta_n(\omega)=\frac{\sum_{l=1}^{L_n^2(\omega)} \log a_{i_lj_lk_l}}{\sum_{l=1}^n \log c_{i_l}}
  \underset{n\to\infty}{\longrightarrow} 1.
\]
Given $t\in (\underline{t}, \overline{t})$ and $\rho\in (0,1]$,
consider the probability vector ${\bf{p}}(t,\rho)$, such that
$t({\bf{p}}(t,\rho))=t$, given by Lemma \ref{lemapij}.
Applying (\ref{vert1}) to ${\bf{p}}(t,\rho)$ we obtain
\begin{align*}
  &d_{{\bf{p}}(t,\rho),n}(\omega)=
  \lambda_1({\bf{p}}(t,\rho))+\beta_n(\omega) \lambda_2({\bf{p}}(t,\rho))
   +\eta_n(\omega) t\label{vert2}\\
  &+\frac{\rho \sum_{l=1}^n \log \gamma_{i_l}(t,\rho)-
  \sum_{l=1}^{L_n^1(\omega)} \log \gamma_{i_l}(t,\rho)}{\sum_{l=1}^n \log c_{i_l}}\notag\\
  &+\frac{\alpha(t,\rho) \sum_{l=1}^{L_n^1(\omega)} \log \Bigl(\sum_k a_{i_lj_lk}^{t}\Bigr)
  -\sum_{l=1}^{L_n^2(\omega)} \log \Bigl(\sum_k a_{i_lj_lk}^{t}\Bigr)}
  {\sum_{l=1}^n \log c_{i_l}}, \notag
\end{align*}
where
\[
    \gamma_i(t,\rho)=\sum_{j=1}^{m_i} b_{ij}^{\lambda_2({\bf{p}}(t,\rho))}
    \Bigl(\sum_{k=1}^{m_{ij}} a_{ijk}^{t}\Bigr)^{\alpha(t,\rho)}.
\]
We choose $\rho=\frac{\log  c}{\log b}$ and $t$ such that $\alpha(t,\rho)=\frac{\log b}{\log a}$.
Using estimates (\ref{est1}), (\ref{est2}) and the fact that we are considering an $(a,b,c;\varepsilon)$-sponge, we get
\begin{align*}
  &d_{{\bf{p}}(t,\rho),n}(\omega)\le
  \lambda_1({\bf{p}}(t,\rho))+\beta_n(\omega) \lambda_2({\bf{p}}(t,\rho))
   +\eta_n(\omega) t\\
  &+\frac{\frac{\log  c}{\log b} \sum_{l=1}^n \log \left(\sum_{j=1}^{m_{i_l}} m_{i_lj}^{\frac{\log b}{\log a}} \right)
  - \sum_{l=1}^{\lfloor n \frac{\log c}{\log b} \rfloor} \log \left(\sum_{j=1}^{m_{i_l}} m_{i_lj}^{\frac{\log b}{\log a}} \right)}{\sum_{l=1}^n \log c_{i_l}}\\
    &+\frac{\frac{\log b}{\log a} \sum_{l=1}^{\lfloor n \frac{\log c}{\log b} \rfloor} \log m_{i_lj_l}
  -\sum_{l=1}^{\lfloor n \frac{\log c}{\log a} \rfloor} \log m_{i_lj_l}} {\sum_{l=1}^n \log c_{i_l}} + \tilde{C}\varepsilon + \frac{D}{n}, \notag
\end{align*}
for some constants $\tilde{C}$ and $D$ (depending only on $a$, $b$, $c$ and $\mathcal{I}$) and for every $\varepsilon>0$ sufficiently small. By \cite[Lemma 4.1]{11} we have that
\begin{align*}
\limsup_{n\to\infty} \frac{1}{n}  \Bigl\{&\frac{\log  c}{\log b} \sum_{l=1}^n \log \left(\sum_{j=1}^{m_{i_l}} m_{i_lj}^{\frac{\log b}{\log a}} \right) - \sum_{l=1}^{\lfloor n \frac{\log c}{\log b} \rfloor} \log \left(\sum_{j=1}^{m_{i_l}} m_{i_lj}^{\frac{\log b}{\log a}} \right) \\
  &+\frac{\log b}{\log a} \sum_{l=1}^{\lfloor n \frac{\log c}{\log b} \rfloor} \log m_{i_lj_l} -\sum_{l=1}^{\lfloor n \frac{\log c}{\log a} \rfloor} \log m_{i_lj_l} \Bigr\} \ge0,
\end{align*}
so
\[
\liminf_{n\to\infty} d_{{\bf{p}}(t,\rho),n}(\omega)\le  \lambda_1({\bf{p}}(t,\rho))+ \lambda_2({\bf{p}}(t,\rho))
   + t + \tilde{C}\varepsilon \le \mathrm{VP}(\Lambda) + \tilde{C}\varepsilon.
\]

Now we deal with hypothesis (\ref{hip}). Given $\eta>0$, since the quantitaties in this lemma depend continuously on the numbers
$a_{ijk}$, we can substitute $a_{ijk}$ by arbitrarily close numbers $\tilde{a}_{ijk}$ that satisfy hypothesis (\ref{hip}) and so that at the end we obtain $\liminf_{n\to\infty} d_{{\bf{p}}(t,\rho),n}(\omega)\le \mathrm{VP}(\Lambda) + \tilde{C}\varepsilon + \eta$. Since $\eta$ can be made arbitrarily close to zero by making $\tilde{a}_{ijk}$ sufficiently close to $a_{ijk}$, we get the desired result.
\end{proof}

\begin{lemma}\label{lemafinal}
Let $0<a\le b\le c\le1$ and assume $\underline{t}<\overline{t}$. There exists a positive constant $\tilde{C}$  (depending only on $a,b,c$ and $\mathcal{I}$) such that, for every $\varepsilon>0$ sufficiently small, if $\Lambda$ is an $(a,b,c; \varepsilon)$-sponge then
\begin{equation*}
   \hd\Lambda\le\mathrm{VP}(\Lambda) + \tilde{C}\varepsilon.
\end{equation*}
\end{lemma}
\begin{proof}
Let $\xi>0$.
Consider the \emph{approximate cubes of order $n$} given by $B_n(z)=\chi(B_n(\omega))$
where $\omega\in \chi^{-1}(z)$, $z\in\Lambda$, $n\in\mathbb{N}$. Then it follows from
Lemma \ref{lemasup} that
\begin{equation}\label{dimp}
  \forall_{z\in\Lambda} \;\forall_{N\in\mathbb{N}}\;
  \exists_{n>N}\;\exists_{{\bf{p}}\in\mathcal{P}} :\,
  \frac{\log \mu_{{\bf{p}}}(B_n(z))} {\log |B_n(z)|}\le \mathrm{VP}(\Lambda)+\tilde{C}\varepsilon+\xi.
\end{equation}
Given $\delta, \eta>0$, we shall build a cover $\mathcal{U}_{\delta,\eta}$ of $\Lambda$ by sets
with diameter $<\eta$ such that
\[
    \sum_{U\in\mathcal{U}_{\delta,\eta}} |U|^{\mathrm{VP}(\Lambda)+\tilde{C}\varepsilon+\xi+2\delta}\le
    \sqrt{3}\, (\max a_{ijk}^{-1})\, M_\delta
\]
where $M_\delta$ is an integer depending on $\delta$ but not on
$\eta$. This implies that $\hd\Lambda\le\mathrm{VP}(\Lambda)+\tilde{C}\varepsilon+\xi+2\delta$ which gives what we want because
$\xi$ and $\delta$ can be taken arbitrarily small. Let $c=\max\,c_{i}<1$. It is
clear that there exists a finite number of Bernoulli measures
$\mu_1,...,\mu_{M_\delta}$ such that
\[
      \forall_{{\bf{p}}}\; \exists_{k\in\{1,...,M_\delta\}} :\,
      \frac{\mu_{{\bf{p}}}(B_n)}{\mu_k(B_n)}\le c^{-\delta n}
\]
for all approximate cubes of order $n$, $B_n$. By (\ref{dimp}),
we can build a cover of $\Lambda$ by approximate cubes $B_{n(z^i)},\,i=1,2,...$ that are disjoint
and have diameters $<\eta$, such that
\[
    \mu_{{\bf{p}}^i}(B_{n(z^i)})\ge |B_{n(z^i)}|^{\mathrm{VP}(\Lambda)+\tilde{C}\varepsilon+\xi+\delta}
\]
for some probabilty vectors ${\bf{p}}^i$. It follows that
\begin{align*}
  \sum_i |B_{n(z^i)}|^{\mathrm{VP}(\Lambda)+\tilde{C}\varepsilon+\xi+2\delta} &\le \sum_i \mu_{{\bf{p}}^i}
   (B_{n(z^i)})\, |B_{n(z^i)}|^\delta\\
  &\le \sum_i \mu_{k_i}(B_{n(z^i)}) \,c^{-\delta n(z^i)}
  \, \sqrt{3}\, (\max a_{ijk}^{-1})\, c^{\delta n(z^i)} \\
  &\le \sqrt{3}\, (\max a_{ijk}^{-1})\,\sum_{k=1}^{M_\delta} \sum_i \mu_k(B_{n(z^i)})
  \le\sqrt{3}\, (\max a_{ijk}^{-1})\, M_\delta
\end{align*}
as we wish.
\end{proof}
The upper estimate of Corollary \ref{cor} follows from Lemma \ref{lemafinal}. As in the beginning of this section, we have that
 $\mathrm{VP}(\Lambda_{a,b,c;\varepsilon})\le \mathrm{VP}(\Lambda_{a,b,c}) + D\varepsilon$ for some positive constant $D$
(depending only on $a,b,c$ and $\mathcal{I}$) and for every $\varepsilon>0$ sufficiently small. This together with Lemma \ref{lemafinal} gives the upper estimate in Theorem \ref{cont}.

\end{document}